\documentclass[12pt,amsfonts]{article}
\usepackage{amsmath,amsthm, amssymb, amstext}
\newcommand {\bc} {\begin{center}}
\newcommand {\ec} {\end{center}}
\begin{document}

\title{Growth of solutions for QG and 2D Euler equations}

\author{Diego Cordoba \\
{\small Department of Mathematics} \\
{\small University of Chicago} \\
{\small 5734 University Av, Il 60637} \\
{\small Telephone: 773 702-9787, e-mail: dcg@math.uchicago.edu} \\
{\small and} \\
Charles Fefferman\thanks{Partially supported by NSF grant DMS 0070692.}\\ 
{\small Princeton University} \\
{\small Fine Hall, Washington Road, NJ 08544} \\
{\small Phone: 609-258 4205, e-mail: cf@math.princeton.edu} \\
}

\date{January 17 2001}
%\address{Princeton University}
\maketitle

\markboth{QG and 2D Euler equations}{D.Cordoba and C.Fefferman}

\newtheorem {Thm}{Theorem}
\newtheorem {Def}{Definition}
\newtheorem {Lm}{Lemma}
\newtheorem {prop}{Proposition}
\newtheorem {Rem}{Remark} 
\newtheorem {Cor}{Corollary}
\def\cal{\mathcal}
\newtheorem {Ack*}{Acknowledgments}

\section{Abstract}
   We study the rate of growth of sharp fronts of the Quasi-geostrophic equation and 2D incompressible Euler equations.. The development of sharp fronts are due to a mechanism that piles up  level sets very fast. Under a semi-uniform collapse, we obtain a lower bound on the minimum distance between the level sets.

\section{Introduction}

 The work of Constantin-Majda-Tabak [1] developed an analogy between the Quasi-geostrophic and 3D Euler equations. Constantin, Majda and Tabak  proposed a candidate for a singularity for the Quasi-geostrophic equation. Their numerics showed evidence of a blow-up for a particular initial data, where the level sets of the temperature contain a hyperbolic saddle. The arms of the saddle tend to close in finite time, producing a a sharp front. Numerics studies done later by Ohikitani-Yamada [8] and Constantin-Nie-Schorgofer [2], with the same initial data, suggested that instead of a singularity the derivatives of the temperature where increasing as double exponential in time. 

  The study of collapse on a curve was first studied in [1] for the Quasi-geostrophic equation where they considered a simplified ansatz for classical frontogenesis with trivial topology. At the time of collapse, the scalar $\theta$ is discontinues across the curve $x_2 = f(x_1)$ with different limiting values for the temperature on each side of the front.  They show that under this topology the directional field remains smooth up to the collapse, which contradicts the following theorem proven in [1]: 
\begin{eqnarray*}
 \text{If locally the direction field remains smooth as t}\ \ \ \ \\
\text{ approaches $T_*$, then no finite singularity is possible}\ \ \ \ \\
 \text{ as t approaches $T_*$.}\ \ \ \ \quad \quad \quad \quad \quad \quad \quad \quad 
\end{eqnarray*} 
The simplified ansatz with trivial topology studied in [1] does not describe a hyperbolic saddle.

Under the definition of a simple hyperbolic saddle, in [3], it was shown that the angle of the saddle can not decrease faster than a double exponential in time. 

The criterion obtained in [5] for a sharp front formation for a general two dimensional incompressible flow is :
\begin{eqnarray*}
\text{A necessary condition to have a sharp front at time T is}\\ 
\int_{0}^{T}|u|_{L^{\infty}}(s) ds = \infty\ \quad \quad \quad \quad \ \ \ \ \quad
\end{eqnarray*}
For the Quasi-geostrophic equation it is not known if the quantity $\int_{0}^{T}|u|_{L^{\infty}}(s) ds$ diverges or not. And the criterion does not say how fast the arms of a saddle can close.

In this paper we do not assume anything on the velocity field, and we show that under a semi-uniform collapse the distance between two level curves cannot decrease faster than a double exponential in time. The semi-uniform collapse assumption greatly weakens the assumptions made in [1] for an ansatz for classical frontogenesis, and the simple hyperbolic saddle in [3].

In the case of 2D incompressible Euler equation we are interested in the large time behavior of solutions.

The two equations we discus in this paper, have in common the property that a scalar function is convected by the flow, which implies that the level curves are transported by the flow. The possible singular scenario is due to level curves approaching each other very fast which will lead to a fast growth on the gradient of the scalar function. Below we study the semi-uniform collapse of two level sets on a curve. By semi-uniform collapse we mean that the distance of the two curves in any point are  comparable.

 The equations we study are as follows: 
 
\bc
{\underline{The Quasi-geostrophic (QG) Equation}} 
\ec
 
Here the unknowns are a scalar $\theta(x,t)$ and a velocity field $u(x,t) = (u_1(x,t), u_2(x,t)) \in R^2$, defined for $t\in[0,T^*)$ with $T^*\leq \infty$, and for $x \in \Omega$ where $\Omega = R^2$ or $R^2/Z^2$.  The equations for $\theta$, u are as follows  
\begin{eqnarray}  
\left (\partial_t + u\cdot\nabla_x \right ) \theta = 0 \\ 
u = \nabla_{x}^{\perp}\psi\ \ and \ \ \psi = (-\triangle_x)^{-\frac{1}{2}}\theta, \nonumber
\end{eqnarray} 
where $\nabla_{x}^{\perp} f = (-\frac{\partial f}{\partial x_2}, \frac{\partial f}{\partial x_1})$ for  scalar functions f. The initial condition is $\theta(x,0) = \theta_0(x)$ for a smooth initial datum $\theta_0$. 

\bc 
{\underline{The  Two-Dimensional Euler Equation}}
\ec 
 
The unknown is an incompressible velocity field u(x,t) as above with vorticity denoted by $\omega$. The 2D Euler equation may be written in the form  
\begin{eqnarray}  
\left (\partial_t + u\cdot\nabla_x \right ) \omega = 0 \\
 u = \nabla_{x}^{\perp}\psi\ \ and \ \ \psi = (-\triangle_x)^{-1}\omega, \nonumber
\end{eqnarray} 
with u(x,0) equal to a given smooth divergence free $u_0(x)$.

\section{Results}

Asssume that q = q(x,t) is a solution to (1) or (2), and that a level curve of q can be parameterized by

\begin{eqnarray}
x_2=\phi_{\rho}(x_1,t)\ \ for\ \ x_1\in[a,b] \label{eq:1}
\end{eqnarray}
with $\phi_{\rho}\in C^1([a,b]\cap [0,T^*))$, in the sense that 
\begin{eqnarray}
q(x_1,\phi_{\rho}(x_1,t), t) = G(\rho)\ \ for\ \ x_1\in[a,b],\label{eq:2}
\end{eqnarray}
and for certain $\rho$ to be specified below.

The stream function $\psi$ satisfies 
\begin{eqnarray}
\nabla^{\perp}\psi=u.
\end{eqnarray} 

From (3) and (4), we have
\begin{eqnarray}
\frac{\partial q}{\partial x_1} + \frac{\partial q}{\partial x_2} \frac{\partial\phi_{\rho}}{\partial x_1} = 0 \label{eq:3}
\end{eqnarray}
\begin{eqnarray}
\frac{\partial q}{\partial t} + \frac{\partial q}{\partial x_2} \frac{\partial\phi_{\rho}}{\partial t} = 0 \label{eq:4}
\end{eqnarray}

By (1), (2), (5), (6) and (7) we obtain
\begin{eqnarray*}
\frac{\partial\phi_{\rho}}{\partial t} & = & -\frac{\frac{\partial q}{\partial t}}{\frac{\partial q}{\partial x_2}} = \frac{<-\frac{\partial\psi}{\partial x_2},\frac{\partial\psi}{\partial x_1}>\cdot <\frac{\partial q}{\partial x_1},\frac{\partial q}{\partial x_2} >}{\frac{\partial q}{\partial x_2}} \\
& = & <-\frac{\partial\psi}{\partial x_2},\frac{\partial\psi}{\partial x_1}>\cdot <\frac{\frac{\partial q}{\partial x_1}}{\frac{\partial q}{\partial x_2}}, 1> \\
& = & <-\frac{\partial\psi}{\partial x_2},\frac{\partial\psi}{\partial x_1}>\cdot <-\frac{\partial\phi_{\rho}}{\partial x_1}, 1>
\end{eqnarray*}
Next 
\begin{eqnarray*}
\frac{\partial}{\partial x_1}\left (\psi(x_1,\phi_{\rho}(x_1,t), t)\right ) & = & \frac{\partial\psi}{\partial x_1} + \frac{\partial\psi}{\partial x_2} \frac{\partial\phi_{\rho}}{\partial x_1} \\
& = & <-\frac{\partial\psi}{\partial x_2},\frac{\partial\psi}{\partial x_1}>\cdot <-\frac{\partial\phi_{\rho}}{\partial x_1}, 1>
\end{eqnarray*}

Therefore
\begin{eqnarray}
\frac{\partial\phi_{\rho}}{\partial t} =  \frac{\partial}{\partial x_1}\left (\psi(x_1,\phi_{\rho}(x_1,t), t)\right ) \label{eq:5}
\end{eqnarray}

With this formula we can write a explicit equation for the change of time of the area between two fixed points a, b and two level curves $(\phi_{\rho_1}, \phi_{\rho_2})$;
\begin{eqnarray}
\frac{d }{d t}\left ( \int_{a}^{b} [\phi_{\rho_2}(x_1,t) -  \phi_{\rho_1}(x_1,t)] dx_1 \right ) \nonumber \\ 
=  \psi(b,\phi_{\rho_2}(b,t), t)  - \psi(a,\phi_{\rho_2}(a,t), t) \nonumber \\
 +  \psi (a,\phi_{\rho_1}(a,t), t) - \psi(b,\phi_{\rho_1}(b,t), t) \label{eq:6}
\end{eqnarray}

 Assume that two level curves  $\phi_{\rho_1}$ and $\phi_{\rho_2}$ collapse when t tends to $T^*$ uniformly in $a\leq x_1\leq b$ i.e.
$$
\phi_{\rho_2}(x_1,t) - \phi_{\rho_1}(x_1,t) \sim \frac{1}{b - a} \int_{a}^{b} [\phi_{\rho_2}(x_1,t) - \phi_{\rho_1}(x_1,t)] dx_1
$$
In other words; the distance between two level sets are comparable for $a \leq x_1 \leq b$. 

Let 
$$
\delta(x_1,t) = |\phi_{\rho_2}(x_1,t) - \phi_{\rho_1}(x_1,t)|
$$
be the thickness of the front.

We define semi-uniform collapse on a curve if (3) and (4) holds and there exists a constant $c$, independent of t, such that
$$
min \delta (x_1,t) \geq c\cdot  max \delta (x_1,t) 
$$
for $a\leq x_1 \leq b$, and for all $t\in [0,T^*)$.

We call the length b-a of the interval [a,b] the length of the front.

Now we can state the following theorem

\begin{Thm} 
For a QG solution with a semi-uniform front, the thickness $\delta(t)$ satisfies 
\begin{eqnarray*} 
\delta(t) > e^{-e^{At + B}} \ \ for\ \ all\ \ t\in[0,T^*). 
\end{eqnarray*}   
Here, the constants A and B may be taken to depend only on the length of the front, the semi-uniformity constant, the initial thickness $\delta(0)$, and the norm of the initial datum $\theta_0(x)$ in $L^1\cap  L^{\infty}$.
\end{Thm}

Proof: From (9) we have

\begin{eqnarray}
|\frac{d}{d t} A(t)| < \frac{C}{b - a}   sup_{a\leq x_1\leq b} |\psi(x_1,\phi_{\rho_2}(x_1,t), t)  - \psi(x_1,\phi_{\rho_2}(x_1,t), t)|
\end{eqnarray}
where
\begin{eqnarray*}
A(t) = \frac{1}{b - a} \int_{a}^{b} [\phi_{\rho_2}(x_1,t) & - &\phi_{\rho_1}(x_1,t)] dx_1,
\end{eqnarray*}
and C is determined by the semi-uniformity constant c.

The estimate of the difference of the value of the stream function at two different points that are close to each other is obtained by writing the stream function as follows;
\begin{eqnarray*}
\psi(x,t) = - \int_{\Omega}\frac{\theta(x + y,t)}{|y|} dy,
\end{eqnarray*}
and this is because $\psi = (-\triangle_x)^{-\frac{1}{2}}\theta$.

Therefore 
\begin{eqnarray*}
  \psi(z_1, t)  - \psi(z_2,t) & = & \int_{\Omega}\theta(y)(\frac{1}{|y - z_1|} - \frac{1}{|y - z_2|}) dy \\ &  = & \int_{|y - z_1| \leq 2\tau} + \int_{2\tau < |y - z_2| \leq k} + \int_{k < |y - z_1| } \\ & \equiv & I_{1} + I_{2} + I_{3}.\end{eqnarray*}
where $\tau = |z_1 - z_2| $.

Furthermore
\begin{eqnarray*}
|I_{1}| & \leq & ||\theta||_{L^{\infty}} \cdot\int_{|y - z_1| \leq 2\tau}(\frac{1}{|y - z_1|} + \frac{1}{|y - z_2|}) dy \\ & \leq & C\tau
\end{eqnarray*}

We define s to be a point in the line between $z_1$ and $z_2$, then $|y - z_1|
\leq 2|y - s|$ and $I_{2}$ can be estimated by
\begin{eqnarray*}
|I_{2}| &\leq& C\tau \cdot \int_{2\tau < |y - z_1| \leq
k}max_{s}|\nabla(\frac{1}{|y - s|})| dy \\ &\leq& C\tau \cdot
\int_{2\tau < |y - z_1| \leq k}max_{s}\frac{1}{|y - s|^{2}} dy \\ &\leq&
C\tau \cdot |\log \tau|
\end{eqnarray*}

We use the conservation of energy to estimate $I_{3}$ by 
\begin{eqnarray*}
|I_{3}| \leq C \cdot \tau
\end{eqnarray*}

Finally, by choosing $\tau = |z_1 - z_2|$ we obtain
\begin{eqnarray}
|\psi(z_1, t)  - \psi(z_2,t)| \leq M|z_1 - z_2||log |z_1 - z_2||
\end{eqnarray} 
where M is a constant that depend on the initial data $\theta_0$. (See details in [3].)

Then we have 
\begin{eqnarray*}
|\frac{d}{d t} A(t)| & \leq & \frac{M}{b - a} sup_{a\leq x_1\leq b}|\phi_{\rho_2}(x_1,t) - \phi_{\rho_1}(x_1,t)||log |\phi_{\rho_2}(x_1,t) - \phi_{\rho_1}(x_1,t)|| \\
& \leq & \frac{C\cdot M}{\cdot(b - a)} |A(t)||log A(t)|
\end{eqnarray*}
and therefore
\begin{eqnarray*}
A(t) >> A(0)e^{-e^{\frac{C\cdot M}{\cdot(b - a)} t}}
\end{eqnarray*}

\begin{Thm} 
For a 2D Euler solution with a semi-uniform front, the thickness $\delta(t)$ satisfies 
\begin{eqnarray*} 
\delta(t) > e^{-[At + B]} \ \ for\ \ all\ \ t\in[0,T^*). 
\end{eqnarray*}  
\end{Thm}

Here, the constants A and B may be taken to depend only on the length of the front, the semi-uniformity constant, the initial thickness $\delta(0)$, and the norm of the initial vorticity in $L^1\cap  L^{\infty}$.

The proof theorem 2 is similar to theorem 1 with the difference that instead of the estimate (11), we have
\begin{eqnarray*}
|\psi(z_1, t)  - \psi(z_2,t)|\leq M|z_1 - z_2|
\end{eqnarray*}
where M is a constant that depend on the initial data $u_0$. (See details in [3].)

Similar estimates can be obtain for 2D ideal  Magneto-hydrodynamics (MHD) Equation, with the extra assumption that $\int_{0}^{T^*}|u|_{L^{\infty}}(s) ds$ is bounded up to the time of the blow-up. This estimates are consequence of applying the Mean value theorem in (10). Nevertheless in the case of MHD these estimates improve the results obtain in [6].

\begin{Ack*} 
 This work was initially supported by the American Institute of Mathematics. 
\end{Ack*}

\end{document}